\documentclass[3p]{elsarticle}

\usepackage{graphicx}
\usepackage{subfigure}
\usepackage{amsfonts}
\usepackage{amssymb}
\usepackage{amsmath}
\usepackage{hyperref}
\usepackage{bm}
\usepackage[bottom]{footmisc}

\newcommand{\R}{\mathbb{R}}
\newcommand{\bo}{\partial}
\newcommand{\closure}[1]{\overline{#1}}

\newcommand{\interior}{\mathrm{int}}

\newcommand{\heading}[1]{\smallskip\par\noindent{\bf #1}}
\def\O{\mathcal{O}{}}

\def\claimqed{\hfill$\diamond$}
\newcommand{\Gsub}[1]{G^{(#1)}}

  \def\calC{{\cal C}}

\newtheorem{theorem}{Theorem}[section]
\newtheorem{lemma}[theorem]{Lemma}
\newtheorem{claim}[theorem]{Claim}

\newtheorem{corollary}[theorem]{Corollary}
\newtheorem{observation}[theorem]{Observation}

\newtheorem{problem}{Problem}
\newproof{proof}{Proof}
\newproof{claim_proof}{Proof of Claim}
\def\int{\hbox{\bf \rm \sffamily INTERVAL}}
\def\kint#1{\hbox{\bf \rm \sffamily $#1$-INTERVAL}}
\def\kunitint#1{\hbox{\bf \rm \sffamily $#1$-UNIT INTERVAL}}

\def\ifa{\hbox{\bf \rm \sffamily INTERVAL FILAMENT}}
\def\str{\hbox{\bf \rm \sffamily STRING}}
\def\outerstr{\hbox{\bf \rm \sffamily OUTER-STRING}}
\def\genusstr#1{\hbox{\bf \rm \sffamily $#1$-GENUS STRING}}
\def\eulergenusstr#1{\hbox{\bf \rm \sffamily $#1$-EULER-GENUS STRING}}

\def\genus#1{\hbox{\bf \rm \sffamily GENUS $#1$}}

\def\boxicity#1{\hbox{\bf \rm \sffamily $#1$-BOX}}
\def\cube#1{\hbox{\bf \rm \sffamily $#1$-CUBE}}
\def\ball#1{\hbox{\bf \rm \sffamily $#1$-BALL}}

\def\grid#1{\hbox{\bf \rm \sffamily $#1$-GRID}}
\def\unitcube#1{\hbox{\bf \rm \sffamily $#1$-UNIT CUBE}}
\def\unitball#1{\hbox{\bf \rm \sffamily $#1$-UNIT BALL}}
\newcommand{\linegr}{\hbox{\bf \rm \sffamily LINE}}

\def\cNP{\hbox{\rm \sffamily NP}}

\def\cPSPACE{\hbox{\rm \sffamily PSPACE}}
\def\cEXPTIME{\hbox{\rm \sffamily EXPTIME}}
\def\cn{\textrm{\rm cn}}
\def\mcn{\textrm{\rm max-cn}}

\newenvironment{packed_itemize}{
	\begin{itemize}
		\setlength{\itemsep}{1pt}
	    \setlength{\parskip}{0pt}
		\setlength{\parsep}{0pt}
}{\end{itemize}}

\begin{document}

\title{Cops and Robbers on Intersection Graphs\tnoteref{conference}}
\tnotetext[conference]{The conference versions of parts of this paper appeared in ISAAC
2013~\cite{ISAAC2013} and ISAAC 2015~\cite{ISAAC2015}. For a structural dynamical diagram of the
results of this paper, see \url{http://pavel.klavik.cz/orgpad/cops\_on\_intersection\_graphs.html}
(supported for Firefox and Google Chrome).  The third, the fourth, and the fifth authors are
supported by CE-ITI (P202/12/G061 of GA\v{C}R), the first, the fourth and the fifth authors are
supported by Charles University as GAUK 196213.}

\author[cunidam]{Tom\'a\v{s} Gaven\v{c}iak}
\ead{gavento@kam.mff.cuni.cz}
\author[lodz]{Przemys{\l}aw Gordinowicz}
\ead{pgordin@p.lodz.pl}
\author[cunicsi]{V\'{\i}t Jel\'{\i}nek}
\ead{jelinek@iuuk.mff.cuni.cz}
\author[cunicsi]{Pavel Klav\'ik}
\ead{klavik@iuuk.mff.cuni.cz}
\author[cunidam]{Jan Kratochv\'il}
\ead{honza@kam.mff.cuni.cz}

\address[cunidam]{Department of Applied Mathematics, Faculty of Mathematics and Physics,\\Charles
		University in Prague, Prague, Czech Republic.}
\address[lodz]{Institute of Mathematics, Technical University of Lodz, {\L}\'{o}d\'{z}, Poland.}
\address[cunicsi]{Computer Science Institute, Charles University in Prague, Prague, Czech Republic.}

\begin{abstract}
The cop number of a graph $G$ is the smallest $k$ such that $k$ cops win the game of cops and robber on $G$.
We investigate the maximum cop number of geometric intersection graphs, which are graphs whose vertices are
represented by geometric shapes and edges by their intersections. We establish the following
dichotomy for previously studied classes of intersection graphs:
\begin{packed_itemize}
\item The intersection graphs of arc-connected sets in the plane (called \emph{string graphs}) have cop number at most 15, and more generally, the intersection graphs of arc-connected subsets of a surface have cop number at most $10g+15$ in case of orientable surface of genus $g$, and at most $10g'+15$ in case of non-orientable surface of Euler genus $g'$. For more restricted classes of intersection graphs, we obtain better bounds: the maximum cop number of interval filament graphs is two, and the maximum cop number of outer-string graphs is between 3 and 4.
\item The intersection graphs of disconnected 2-dimensional sets or
of 3-dimensional sets have unbounded cop number even in very restricted settings.
For instance, we show that the cop number is unbounded on intersection graphs of
two-element subsets of a line, as well as on intersection graphs of 3-dimensional
unit balls, of 3-dimensional unit cubes or of 3-dimensional axis-aligned unit segments.
\end{packed_itemize}
\end{abstract}

\begin{keyword}
intersection graphs\sep
string graphs\sep
graphs on surfaces\sep
interval filament graphs\sep
cop and robber\sep
pursuit games\sep
games on graphs
\end{keyword}

\maketitle

\section{Introduction}\label{sec:introduction}

The game of cops and robber on graphs has been introduced independently by
Quilliot~\cite{quilliot78,quilliot83} and by Winkler and Nowakowski~\cite{WN}. In this paper, we
investigate the game on geometric intersection graphs.

\heading{Rules of the Game.}
The first player, called \emph{the cops}, places $k$ cops on vertices of a graph $G$. Then the
second player, called \emph{the robber}, places the robber on a vertex. Then the players
alternate. In the cops' move, every cop either stays in its vertex, or moves to one of
its neighbors.  More cops may occupy the same vertex. In the robber's move, the robber either stays
in its vertex, or moves to a neighboring vertex.  The game ends when the robber is \emph{captured}
which happens when a cop occupies the same vertex as the robber. The cops win if they are able to
capture the robber. The robber wins if he is able to escape indefinitely.

\heading{Maximum Cop Number.}
For a graph $G$, its \emph{cop number} $\cn(G)$ is the least number $k$ such that $k$ cops have a
winning strategy on $G$. For a class of graphs $\calC$, the \emph{maximum cop number} $\mcn(\calC)$
is the maximum cop number $\cn(G)$ of a connected graph $G \in \calC$, possibly $+\infty$.  The
restriction to connected graphs is standard: if $G$ has connected components $C_1,\dots,C_k$, then
$\cn(G) = \sum_{i=1}^k \cn(C_i)$. Therefore, a graph class closed under disjoint union cannot have a
bounded maximum cop number if we omit this restriction. Throughout the paper, we only work with
connected graphs.

\heading{Known Results.}
Graphs of the cop number one were characterized already by Quilliot~\cite{quilliot83} and by
Nowakowski and Winkler~\cite{WN}.  These are the graphs whose vertices can be linearly ordered $v_1,
v_2, \ldots , v_n$ so that each $v_i$ for $i\geq 2$ is a corner of $G[v_1,\ldots,v_i]$, i.e., $v_i$
has a neighbor $v_j$ for some $j<i$ such that $v_j$ is adjacent to all other neighbors of $v_i$.
Andreae~\cite{bounded_degree} proved that $k$-regular graphs have the maximum cop number equal
$+\infty$ for all $k \ge 3$.

For $k$ part of the input, deciding whether the cop number of a graph is at most $k$ has been shown
to be \cNP-hard~\cite{FGKNS}, \cPSPACE-hard~\cite{Mamino13} and very recently
\cEXPTIME-complete~\cite{Kinnersley15}, confirming a 20 years old conjecture of Goldstein and
Reingold~\cite{GoldsteinR95}.  In order to test whether $k$ cops suffice to capture the robber on an
$n$-vertex graph, we can search the game graph which has $\O(n^{k+1})$ vertices to find a winning
strategy for cops. In particular, if $k$ is a fixed constant, this algorithm runs in polynomial
time.

For general graphs on $n$ vertices, it is known that at least $\sqrt{n}$ cops may be needed (e.g.,
for the incidence graph of a finite projective plane~\cite{pralat}). Meyniel conjecture states that the cop
number of a connected $n$-vertex graph is $\O(\sqrt{n})$. For more details and results, see the
book~\cite{Now}.

\heading{Geometrically Represented Graphs.}
We want to argue that the geometry of a graph class heavily influences the maximum cop number. For
instance, the classical result of Aigner and Fromme~\cite{aigner_fromme} shows that the maximum cop
number of planar graphs is 3. This result was generalized to graphs of bounded genus by
Quilliot~\cite{Q2} and improved by Schroeder~\cite{Sch}:
\begin{equation} \label{eq:cn_genus}
\cn(G) \le \frac{3}{2}g + 3,
\end{equation}
where $g$ is the (orientable) genus of $G$. For non-orientable surfaces, a similar result was obtained by Clarke et al.~\cite{Clarke2014}: 
\begin{equation} \label{eq:cn_euler_genus}
\cn(G) \le \frac{3}{2}g' + \frac{3}{2},
\end{equation}
where $g'$ is the Euler genus of $G$ (also called the \emph{crosscap number} of the surface $G$ is
drawn on). However, the exact value of the maximum cop number is not known already for toroidal
graphs ($g = 1$). 

We study \emph{intersection representations} in which a graph $G$ is represented by a map $\varphi\colon
V \to 2^X$ for some ground set $X$ such that the edges of $G$ are described by the intersections:
$uv\in E \iff \varphi(u)\cap\varphi(v)\neq\emptyset$. The ground set $X$ and the images of $\varphi$
are usually somehow restricted to get particular classes of intersection graphs. For example, the
well-known interval graphs have $X=\R$ and every $\varphi(v)$ is a closed interval.

All these graph classes admit large cliques, so their genus is unbounded and the
bound~(\ref{eq:cn_genus}) of the maximum cop number does not apply. On the other hand, existence of
large cliques does not imply big maximum cop number since only one cop is enough to guard a maximal clique.
For instance, chordal graphs, which are intersection graphs of subtrees of a tree, may have
arbitrary large cliques but their maximum cop number is 1.

\begin{figure}[!p]
\centering
\includegraphics{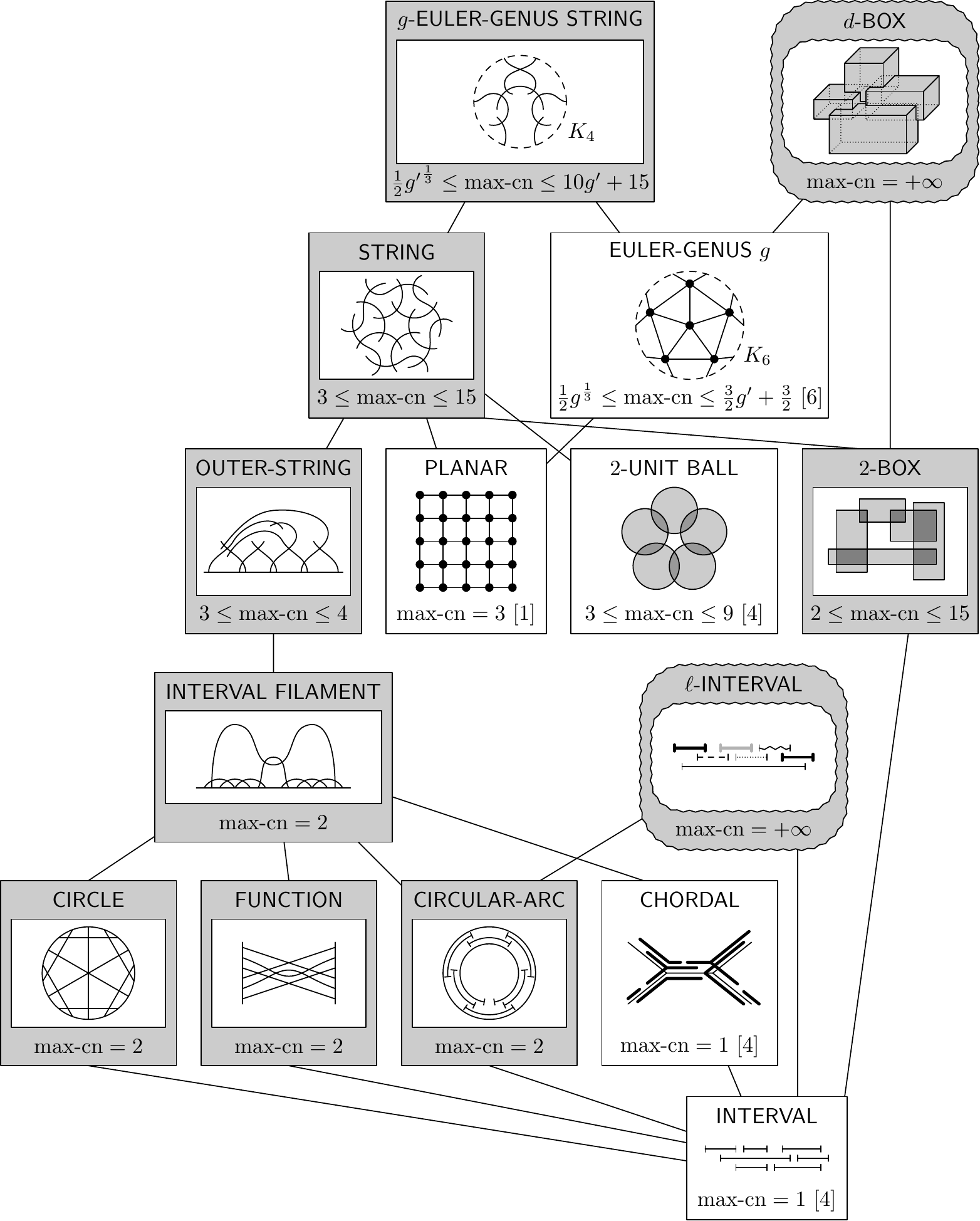}
\caption{The Hasse diagram of inclusions of the considered classes of graphs, together with bounds
on the maximum cop number. The classes with previously known bounds are depicted in white, and the
classes with the bounds proved in this paper are depicted in gray. Bounded boxicity of bounded genus
graphs has been shown in~\cite{EJ13}.}
\label{fig:classes}
\end{figure}

\heading{String Graphs.}
The class of \emph{string graphs} (\str) is the class of intersection graphs of \emph{strings}:
$X=\R^2$ and every $\varphi(v)$ is a bounded curve, i.e., a continuous image of the interval $[0,1]$
in $\R^2$.  It is known that every intersection graph of {\em arc-connected sets} in the plane is a
string graph.  For instance, \emph{boxicity $d$ graphs} (\boxicity{d}), which are intersection
graphs of $d$-dimensional intervals in $\R^d$, are string graphs when $d \le 2$.

The class of \emph{outer-string graphs} (\outerstr) consists of all string graphs having string
representations with each string in the upper half-plane, intersecting the $x$-axis in exactly one
point, which is an endpoint of this string.

The class of \emph{interval filament graphs} (\ifa), introduced by Gavril~\cite{gavril}, consists of
intersection graphs of interval filaments, where an \emph{interval filament} $\varphi(v)$ defined on
an interval $[a,b]$ is the graph of a continuous function $f_v\colon [a,b] \to \R$ such that
$f_v(a)=f_v(b) = 0$ and $f_v(x) > 0$ for all $x \in (a,b)$. In our description, we identify the
filament $\varphi(v)$ with the function $f_v$, i.e., we use $\varphi(v)(x)$ instead of $f_v(x)$.
It holds that $\ifa \subsetneq \outerstr \subsetneq \str$. (The first
inequality is strict from Theorem~\ref{thm:bounded_max_cn}(i) and (ii), the second
follows from~\cite{unbounded_outerstring}.)

Let $\mathbf S$ be an arbitrary orientable surface of genus $g$. We consider a generalization of
string graphs for $X = \mathbf S$ and every $\varphi(v)$ is a bounded curve in $\mathbf S$, and we
denote this class by \genusstr{g}. It can be seen that any intersection graph of arc-connected
subsets of $\mathbf S$ belongs to~\genusstr{g}.  It is known that every graph embeddable to a
surface of genus $g$ can be represented by a contact representation of disks on a suitable Riemann
surface of genus $g$; so $\genus{g} \subsetneq \genusstr{g}$. (It is strict since arbitrarily large
complete graphs belong to $\genusstr{g}$, but not to $\genus{g}$.) Similarly. let \eulergenusstr{g'}
denote intersection graphs of bounded curves on a (possibly non-orientable) surface of Euler genus
$g'$ (equal to the crosscap number of the surface).

Note that while we could work only with Euler genus for both orientable and non-orientable surfaces,
the bounds obtained for orientable genus are significantly lower that for Euler genus (since
$g' \le 2g$ is the best possible general bound).

\heading{Intersection Graphs of Disconnected and Higher Dimensional Sets.} As stated, all intersection
classes of graphs of arc-connected sets in the plane are subclasses of string graphs. Other classes
of intersection graphs are obtained either for disconnected sets, or for sets in dimensions higher
than two.

For a graph $G$, it \emph{line-graph}, denoted by $L(G)$, is the intersection graph of the edges
of~$G$.  Let $\linegr$ denote the class of all line-graphs.  Observe that each line-graph can be
represented as an intersection graph of two-element subsets of a line.  Thus, line-graphs provide a
simple example of intersection graphs of disconnected sets.  As shown by Dudek et al.~\cite{edges},
the cop number of $L(G)$ is related to the cop number of $G$ via the inequalities 
$$\left\lceil\frac{\cn(G)}{2}\right\rceil\le \cn(L(G))\le \cn(G)+1.$$
In particular, the cop number of line-graphs is unbounded.

Many other geometric intersection classes can be seen as generalizations of line-graphs.
Among the most studied are the \emph{$\ell$-interval graphs} (\kint{\ell}) where $X =
\R$ and every $\varphi(v)$ is a union of $\ell$ closed intervals. And for its subclass
\emph{$\ell$-unit interval graphs} (\kunitint{\ell}), all $\ell$ intervals of $\varphi(v)$ have the
length one. It follows that the cop number is unbounded on all these classes.

For sets in higher dimensions, notice that every graph has a representation by 3-dimensional
strings. Therefore, to get interesting classes of graphs, we have to further restrict geometry of
the sets. Aside already described \boxicity{d}, we consider the following classes for $X = \R^d$:
the intersection graphs of axis parallel segments \grid{d}, the intersection graphs of
$d$-dimensional unit cubes \unitcube{d}, the intersection graphs of
$d$-dimensional balls \ball{d}, and the intersection graphs of $d$-dimensional unit balls
\unitball{d}.

\heading{Our Results.}
It has been asked at several occasions, last during the Banff Workshop on Graph Searching in October
2012, whether intersection-defined graph classes (other than interval graphs) have bounded maximum
cop numbers. The classes in question have included circle graphs, intersection graphs of disks in
the plane, graphs of boxicity 2, and others. A recent paper~\cite{unit_disk} shows that the maximum
cop number of intersection graphs of unit disks (\unitball{2}) is between 3 and 9. We solve this
question in a general way by proving a dichotomy for previously studied classes of geometric
intersection graphs in Theorems~\ref{thm:bounded_max_cn} and \ref{thm:unbounded_max_cn}.  For an
overview of the results presented in this paper, see Fig.~\ref{fig:classes}. 

\begin{theorem} \label{thm:bounded_max_cn}
The following bounds for the maximum cop number hold:
\begin{enumerate}
\item[(i)] $\mcn(\ifa) = 2$. 
\item[(ii)] $3 \le \mcn(\outerstr) \le 4$.
\item[(iii)] $3 \le \mcn(\str) \le 15$.
\item[(iv)] ${1 \over 2} g^{1 \over 3} \le \mcn(\genusstr{g}) \le 10g+15$.
\item[(v)] ${1 \over 2} {g'}^{1 \over 3} \le \mcn(\eulergenusstr{g'}) \le 10g'+15$.
\end{enumerate}
\end{theorem}

We note that the strategies of cops in all upper bounds are geometric and their description is
constructive, using an intersection representation of $G$.  If only the graph $G$ is given, we
cannot generally construct these representations efficiently since recognition is \cNP-complete for
string graphs~\cite{kratochvil_string_nphard} and interval filament graphs~\cite{pergel}, and open
for the other classes.  Nevertheless, since the state space of the game has $\O(n^{k+1})$ states and
the number of cops $k$ is bounded by a constant, we can use the standard exhaustive game space
searching algorithm to obtain the following:

\begin{corollary}
There are polynomial-time algorithms computing the cop number and an optimal strategy for the cops
for any interval filament graph in time $\O(n^3)$, outer-string graph in time $\O(n^5)$, string
graph in time $\O(n^{16})$ and a string graph on a surface of a fixed genus $g$ (resp. Euler genus $g'$) in time $\O(n^{10g+16})$ (resp. $\O(n^{10g'+16})$), even when representations are not given.\qed
\end{corollary}

Furthermore, our results can be used as a polynomial-time heuristic to prove that a given graph $G$
is not, say, a string graph, by showing that $\cn(G) > 15$.  For instance, a graph $G$ of girth 5
and the minimum degree at least 16 is not a string graph since $\cn(G) > 15$: in any position of 15
cops with the robber on $v$, at least one neighbor of $v$ is non-adjacent to the cops.

On the other hand, when sets are not arc-connected, we prove that even in very restricted
geometric settings that cop numbers are unbounded. The main lemma states that when we subdivide all
edges of a graph $G$ by a same number of vertices, $\cn(G)$ increases by at most one. Since all
these classes contain certain subdivisions of all graphs or all cubic graphs, we get the following:

\begin{theorem} \label{thm:unbounded_max_cn}
The classes \linegr, \kint{2}, \kunitint{2}, \grid{3}, \boxicity{3}, \unitcube{3}, \ball{3}, and \unitball{3}
have the maximum cop number equal $+\infty$.
\end{theorem}

\heading{Outline.} In Section~\ref{sec:ifa}, we show that $\mcn(\ifa)$ is 2. In
Section~\ref{sec:outerstring}, we show that $\mcn(\outerstr)$ is between 3 and 4. Aigner and
Fromme~\cite{aigner_fromme} show the classical result that one cop can guard a shortest path in any
graph. In Section~\ref{sec:guarding_path}, we extend this result to show that five cops can guard a
shortest path together with its neighborhood.  This is used in Section~\ref{sec:string} to show that
$\mcn(\str)$ is at most 15. In Section~\ref{sec:bounded_genus}, we combine
the previous result with the approach of Quilliot~\cite{Q2} to simultaneously show the bounds for bounded-genus orientable surfaces and bounded-Euler-genus non-orientable surfaces. In Section~\ref{sec:unbounded_max_cn}, we prove that cop numbers are unbounded for
intersection graphs of disconnected or 3-dimensional sets.

\heading{Preliminaries.}
Let $G = (V,E)$ be a graph. For $D \subseteq V$, we let $G[D]$ denote the subgraph of $G$ induced by $D$, and $G -
v = G[V \setminus \{v\}]$. For a vertex $v$, we use the \emph{open neighborhood} $N(v) = \{u : uv
\in E\}$ and the \emph{closed neighborhood} $N[v] = N(v) \cup \{v\}$. Similarly for $V' \subseteq
V$, we put $N[V']= \bigcup_{v \in V'}N[v]$ and $N(V') = N[V'] \setminus V'$.

Let $\varphi\colon V\to 2^{\R^2}$ be a string representation of $G$. Without loss of generality, we may
assume that there is only a finite number of string intersections in the representation, that
strings never only touch without either also crossing each other, or at least one of them ending,
that no three or more strings meet at the same point, and that no string self-intersects.  This
follows from the fact that strings can be replaced by piece-wise linear curves with finite 
numbers of linear segments without affecting their intersection graph. For more details 
see~\cite{KGK}. We always assume and maintain these properties.

Suppose that we have some strategy for the cops.  For a vertex $v \in V$, the robber \emph{cannot
safely move to $v$} if the strategy ensures that he is immediately captured after moving to $v$. Let
$D,P \subseteq V$.  We say that the strategy \emph{guards $P$} if it ensures that the robber cannot
safely move to any vertex in $P$.  We say that the robber is \emph{confined to $D$}, if the strategy
ensures that the robber is immediately captured by moving to any vertex in $V \setminus D$.  Notice
that the robber is confined to $D$ if and only if he stands in $D$ and $N(D)$ is guarded.

\section{Capturing Robber in Interval Filament Graphs} \label{sec:ifa}

In this section, we show that the maximum cop number of interval filament graphs is equal to two, thus
establishing Theorem~\ref{thm:bounded_max_cn}(i).

If a filament $\varphi(u)$ is defined on $[a,b]$, we call $a$ the \emph{left endpoint} and $b$ the
\emph{right endpoint} of $\varphi(u)$. We assume that the filaments have pairwise distinct endpoints
and the defining intervals are always non-trivial.  In the description, we move the cops on the
representation $\varphi$, and we say that a cop takes a filament $\varphi(u)$ if it is placed on the
vertex $u$ which this filament represents. We shall assume that the robber never moves into the
neighborhood of a vertex taken by a cop, and a cop capture the robber immediately if he stands on a
neighboring vertex.

\heading{Filaments and Regions.}
It is important that each filament splits the half-plane into two regions: the unbounded \emph{top
region} and the \emph{bottom region}. A filament $\varphi(u)$ is \emph{nested in} a filament
$\varphi(v)$ if $\varphi(u)$ is contained in the bottom region of $\varphi(v)$. We say that the
robber \emph{is/stays in a region} if he is/stays on filaments entirely contained in this region.
The robber is \emph{confined by} $\varphi(u)$ if a cop is placed in $\varphi(u)$ and the robber is
in the bottom region of $\varphi(u)$.

\begin{lemma} \label{lem:confinment}
Suppose that the robber is confined by $\varphi(u)$. Then he stays in the bottom region of
$\varphi(u)$ as long as there is a cop on $\varphi(u)$.
\end{lemma}

\begin{proof}
To move from one region to another, the robber has to use a filament $\varphi(v)$ which crosses
$\varphi(u)$. But then $v$ is a neighbor of $u$, and the cop captures the robber in the next
turn.\qed
\end{proof}

A filament $\varphi(u)$ is called \emph{top in $x$} if it maximizes the value $\varphi(v)(x)$ over
all filaments $\varphi(v)$ defined for $x$. Suppose that $\ell$ is the left-most and $r$ is the
right-most endpoint of the representation. We have a \emph{sequence of top filaments}
$\bigl\{\varphi(t_i)\bigr\}_{i=1}^k$ as we traverse from $\ell$ to $r$. We note that one filament can appear
several times in this sequence. See Fig.~\ref{fig:top_filaments} for an example.
    
\begin{figure}[t!]
\centering
\includegraphics{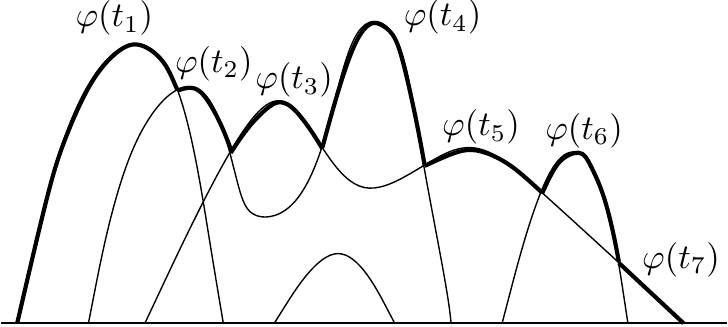}
\caption{An example of a sequence of top filaments. Only the top part of each $\varphi(t_i)$ is
depicted in bold.}
\label{fig:top_filaments}
\end{figure}

Let $\varphi(t_i)$ be top in $x$. Each filament $\varphi(t_i)$ together with the upward ray starting
at $\bigl(x,\varphi(t_i)(x)\bigr)$ separates the half-plane into three regions: the \emph{left
region}, the \emph{bottom region} and the \emph{right region}.  The key property is that there is no
filament intersecting the left and right regions and avoiding $\varphi(t_i)$. Also note that the
division of filaments into the regions is the same for all $x$ in the same top part of
$\varphi(t_i)$.

\begin{lemma} \label{lem:hunting_lemma}
Suppose that a cop stands on $\varphi(t_i)$ and the robber is in the right region of $\varphi(t_i)$. If the cop
moves to a neighboring filament $\varphi(t_j)$ with~$j>i$, the robber cannot move to the
left region of $\varphi(t_j)$.
\end{lemma}

\begin{proof}
Let $\varphi(t_i)$ be on $[a,b]$ and $\varphi(t_j)$ on $[c,d]$ and we have $c < b$.
Suppose that the cop moves from $\varphi(t_i)$ to $\varphi(t_j)$ and the robber
stands on a filament $\varphi(u)$ defined on $[e,f]$. We know that $c < b < e < f$, so $\varphi(u)$ does
not intersect the left region of $\varphi(t_j)$. And since $\varphi(t_j)$ is top, there is no path
going to the left region which avoids $\varphi(t_j)$. So the robber cannot move there.\qed
\end{proof}

\heading{Proof of Theorem~\ref{thm:bounded_max_cn}(i).}
We are ready to prove that the maximum cop number of interval filament graphs is equal to two.

\begin{proof}[Theorem~\ref{thm:bounded_max_cn}({\rm i})]
Since interval filament graphs contain all cycles $C_n$, one cop has no winning strategy. (Note that
$C_4$ is a circle, circular-arc and function graph.) Therefore two cops are necessary.

We describe a strategy how to capture a robber with two cops. We call one cop the \emph{guard}, and
the other one the \emph{hunter}. The strategy proceeds in phases. Every phase starts with both cops
on a filament $\varphi(u)$ such that the robber is confined by it. The guard stays on $\varphi(u)$
till the robber is either captured, or confined by the hunter in some filament $\varphi(v)$ nested in
$\varphi(u)$. By Lemma~\ref{lem:confinment}, the robber can only move in the bottom region
of $\varphi(u)$. If the confinement by $\varphi(v)$ happens, then the guard moves to the filament
$\varphi(v)$ taken by the hunter, ending the phase.  In the next phase the hunter proceeds with
capture the robber inside the bottom region of $\varphi(v)$.

For the initial phase, we imagine that the guard takes some imaginary filament above all filaments
of $\varphi$ so the robber is confined to its bottom region, i.e., to the entire graph $G$.
We can choose both cops to start the first phase at a filament $\varphi(v)$ with left-most left endpoint
and therefore top in $G$.

Suppose that we are in a phase where the guard is placed on $\varphi(u)$.  Let $G_u$ be the
subgraph of $G$ induced by the vertices whose filaments are nested in $\varphi(u)$, and let $C_u$ be
the connected component of $G_u$ containing the vertex occupied by the robber. Since the guard stays
at $\varphi(u)$ till the robber is confined in some nested $\varphi(v)$, the strategy ensures that
the robber must remain in $C_u$ throughout the phase, because any vertex in $N(C_u)$ is adjacent to
the vertex~$u$ guarded by the guard.

Let $\bigl\{\varphi(t_i)\bigr\}_{i=1}^k$ be the sequence of top intervals in the restriction of $\varphi$ to
the vertices of~$C_u$.  The hunter first goes to $\varphi(t_1)$. When he arrives to $\varphi(t_1)$,
the robber cannot be in the left region of $\varphi(t_1)$ since there is no filament of $C_u$
contained there. Now suppose that the hunter is in $\varphi(t_i)$ and assume the induction
hypothesis that the robber is not in the left region of $\varphi(t_i)$. If the robber is confined in
$\varphi(t_i)$, the phase ends with the guard moving towards $\varphi(t_i)$. If the robber is in the
right region of $\varphi(t_i)$, the hunter moves to the neighbor $\varphi(t_j)$ with maximal index
$j$. By Lemma~\ref{lem:hunting_lemma} the robber cannot move to the left region of $\varphi(t_j)$ so
he is either in the bottom, or the right region. The robber cannot stay in the right regions forever
since $\varphi(t_k)$ has no filament of $C_u$ contained in the right region, so eventually the
robber is confined in $\varphi(t_i)$ or captured directly.

Since there are only finitely many filaments nested in each other, the strategy proceeds in finitely
many phases and the robber is captured.\qed
\end{proof}

With a small modification, we can prove that this strategy captures the robber in $\O(n)$ turns.
Suppose that initially both cops are placed in the filament with the left-most endpoint $\ell$ and
there are $p$ phases.  Let $C_i$ be the graph the robber is confined to by the guard on $u_i$ in
the phase $i$, so $C_1=G$ and let $C_{p+1}=\emptyset$. Let $D_i=C_i\setminus C_{i+1}$, and note that
$D_i$ contains all top filaments of $C_i$.

During the $i$-th phase the hunter moves to any top filament of $C_i$ in at most 2 moves (note that
there must be a filament in $G$ which simultaneously intersects $\varphi(u_i)$ and a top filament of
$C_i$), then to the left-most top filament of $C_i$ in at most $|D_i|$ moves using a shortest path
in $D_i$, then takes at most $|D_i|$ steps over the top filaments of $C_i$.  Finally, when the
hunter confines the robber in $C_{i+1}$, it takes the guard at most $|D_i|+2$ steps to get to
$u_{i+1}$ by a similar argument. Since $\sum |D_i|=n$ and the number of phases is also bounded by
$n$, we have used $\O(n)$ turns.

\section{Capturing Robber in Outer-String Graphs} \label{sec:outerstring}

In this section, we prove that the maximum cop number of outer-string graphs is between 3 and 4,
thus establishing Theorem~\ref{thm:bounded_max_cn}(ii). Our strategy is similar to the one described
in Section~\ref{sec:ifa}.

\heading{String Pairs and Regions.}
For a given outer-string representation of $G$, let $v_1,\dotsc,v_n$ be the ordering of the vertices
of $G$ by the $x$-coordinates of the unique intersection of $\varphi(v_i)$ with the $x$-axis. We
say that $v_i$ is \emph{on the left} of $v_j$ and $v_j$ is \emph{on the right} of $v_i$ if $i<j$.

Every pair of intersecting outer-strings $(v_i,v_j)$ divides the half-plane into at least two
regions: the unbounded \emph{top region}, the \emph{bottom region} incident with an interval of the
$x$-axis, and possibly several \emph{middle regions}. The middle regions do not play any role in our
strategy since no string is entirely contained in them. The strings entirely contained in the bottom
region are \emph{surrounded by $\varphi(v_i)$ and $\varphi(v_j)$}, a robber on a vertex surrounded
by $\varphi(v_i)$ and $\varphi(v_j)$, each occupied by a cop, is \emph{confined by $\varphi(v_i)$
and $\varphi(v_j)$}.

The following lemma can be proved the same way as Lemma~\ref{lem:confinment}:

\begin{lemma} \label{lem:outerstring_confinment}
Suppose that the robber is confined by $\varphi(v_i)$ and $\varphi(v_j)$.  Then he stays in the
bottom region of $(v_i,v_j)$ as long as there are cops on $v_i$ and $v_j$.\qed
\end{lemma}

The strings partition the upper half-plane into several regions, of which exactly one is unbounded.
We say that a string $\varphi(x)$ is \emph{external}, if it has at least one point on the boundary
of the unbounded region.  We have a \emph{sequence of external strings}
$\bigl\{\varphi(x_i)\bigr\}_{i=1}^k$ sorted by their appearance on the boundary of the unbounded
region, each external string may appear several times in the sequence.  See
Fig.~\ref{fig:external_strings} for an example.

\begin{figure}[t!]
\centering
\includegraphics{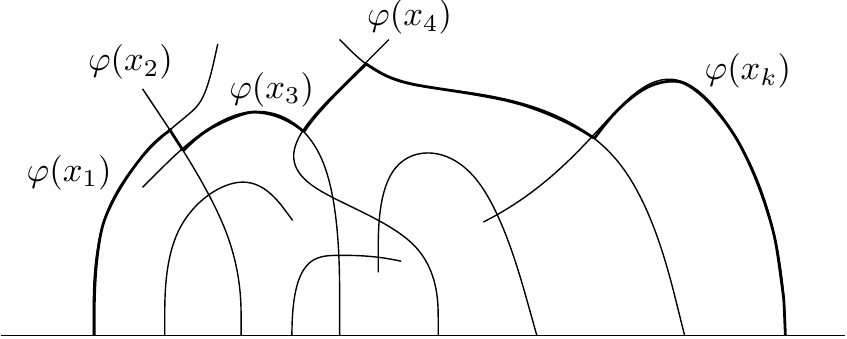}
\caption{An outer-string representation with the sequence of external strings depicted in bold.}
\label{fig:external_strings}
\end{figure}

\begin{lemma} \label{lem:outerstring_external}
Suppose that a cop is placed in an external string $\varphi(x_i)$ and the robber is in some non-intersecting string on
the right of $\varphi(x_i)$. If the cop stays in $\varphi(x_i)$, the robber cannot move to a string
on the left of $\varphi(x_i)$.
\end{lemma}

\begin{proof}
Observe that $\varphi(x_i)$ separates non-intersecting strings on the left of it from those on the
right. Thus, to get to a string on the left of $\varphi(x_i)$, the robber has to move to $N[x_i]$ and the cop captures him.\qed
\end{proof}

\heading{Proof of Theorem~\ref{thm:bounded_max_cn}(ii).}
We are ready to prove that the maximum cop number of outer-string graphs is equal to three or four.

\begin{proof}[Theorem~\ref{thm:bounded_max_cn}({\rm ii})]
Figure~\ref{fig:mrizka} shows a connected outer-string graph with the cop number 3.  It remains to
show that four cops are always sufficient.

\begin{figure}[b!]
\hfil\includegraphics[width=0.3\textwidth]{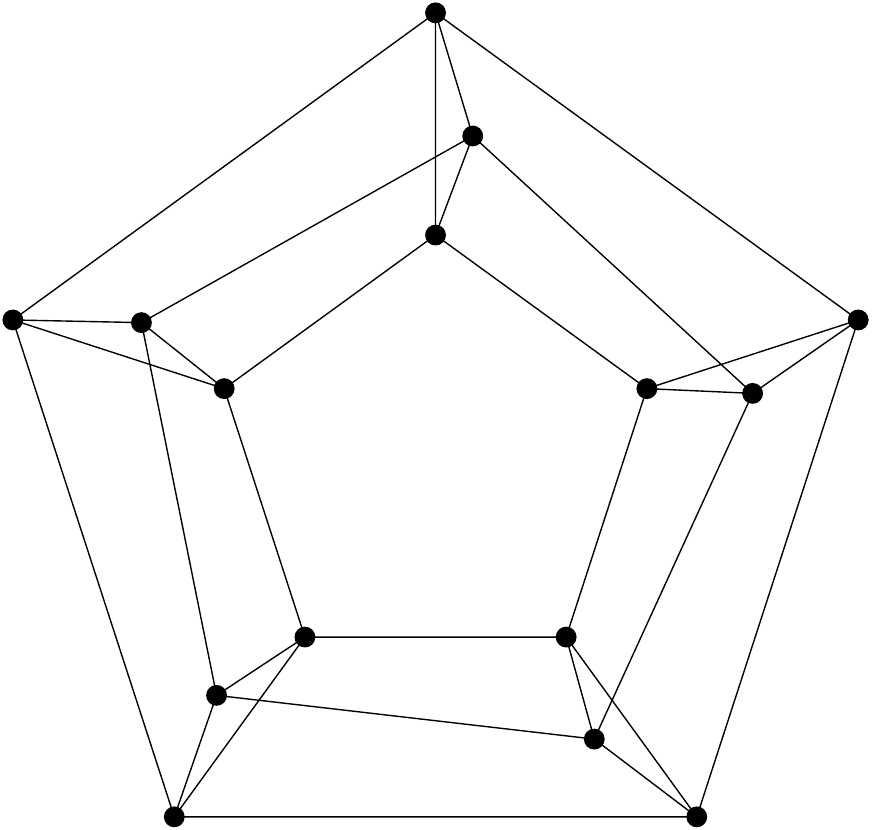}\hfil
\includegraphics[width=0.6\textwidth]{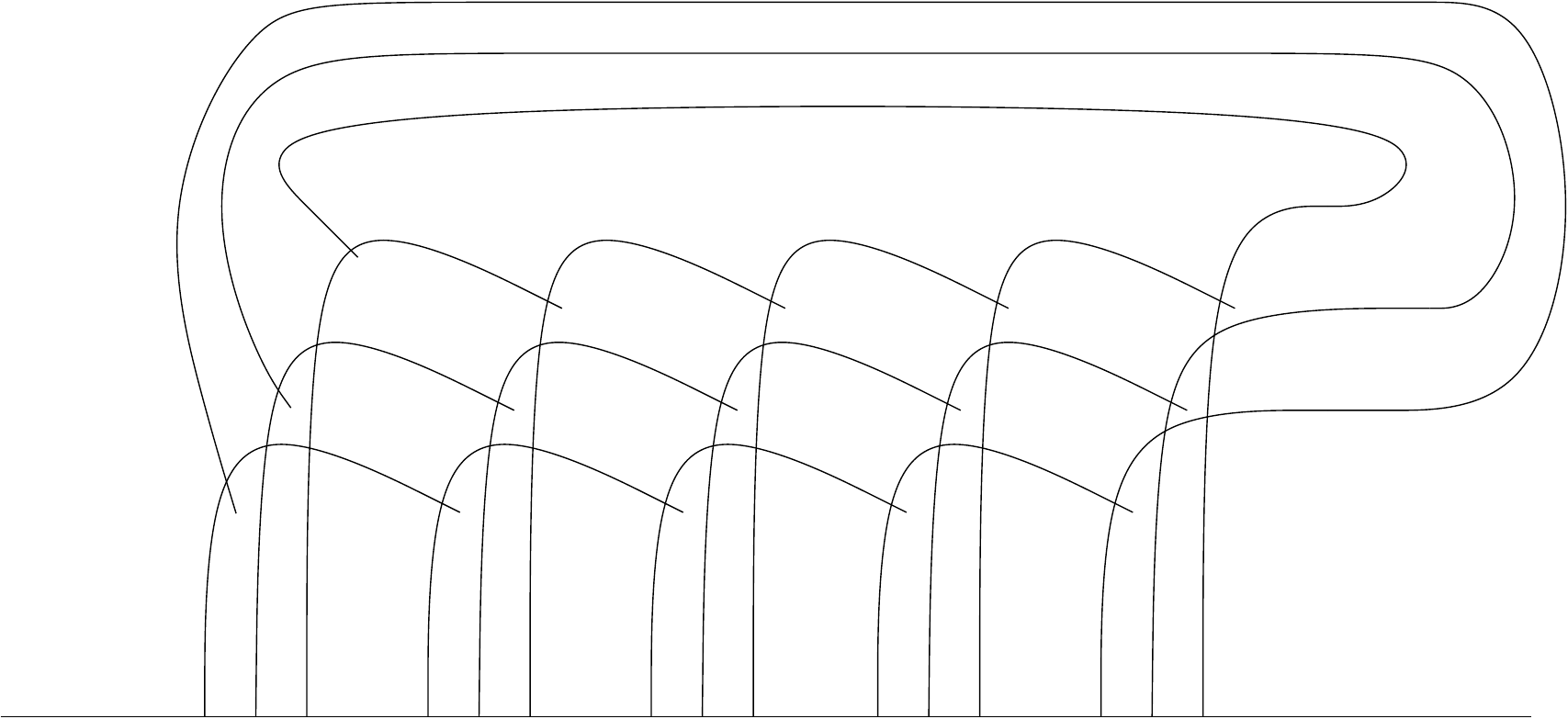}\hfil
\caption{The 3-by-5 toroidal grid and its outer-string representation.  Its cop number is three
since in any position of two cops with at least one cop adjacent to the robber, there is at least
one safe vertex adjacent to the robber.}
\label{fig:mrizka}
\end{figure}

Among the four cops, there are two \emph{guards} and two \emph{hunters}. The strategy is divided
into phases. During each phase, the two guards stand on a pair of intersecting strings
$\varphi(v_\ell)$ and $\varphi(v_r)$ confining the robber.
For the initial phase, we imagine that the guards take some imaginary strings around the entire
representation, so the robber is confined to their bottom region, i.e., to the entire graph $G$.
We can choose all cops to start the first phase in the leftmost string $\varphi(v_1)$ which is
external.

In each phase, let $\varphi(v_\ell)$ and $\varphi(v_r)$ be the pair of adjacent strings occupied by
the guards, confining the robber.  Let $G_{\ell,r}$ be the subgraph induced by the strings entirely
contained in the bottom region, and let $C_{\ell,r}$ be the connected component of $G_{\ell,r}$
containing the vertex with the robber. By Lemma~\ref{lem:outerstring_confinment}, the robber is
confined to $C_{\ell,r}$.  The hunters move in $C_{\ell,r}$ and either capture the robber, or make
the robber confined to a smaller subgraph by taking two intersecting strings $\varphi(v_{\ell'})$
and $\varphi(v_{r'})$ such that the robber is placed in their bottom region.  Then the guards move
to $\varphi(v_{\ell'})$ and $\varphi(v_{r'})$, and the next phase begins.

Let $\bigl\{\varphi(x_i)\bigr\}_{i=1}^k$ be the sequence of external strings in the restriction of
$\varphi$ to the vertices of~$C_{\ell,r}$. The hunters start by taking $\varphi(x_1)$ and its
rightmost neighbor $\varphi(x_i)$.  Suppose that after several moves the hunters take intersecting
external strings $\varphi(x_p)$ and $\varphi(x_q)$.  We want the strategy to preserve the
following property: either the robber occupies a string between $\varphi(x_p)$ and $\varphi(x_q)$,
or he occupies a string on the right of $\varphi(x_q)$.  It is certainly satisfied in the beginning
since there are no strings on the left of $\varphi(x_1)$.

Suppose that the property holds when the hunters occupy $\varphi(x_p)$ and $\varphi(x_q)$, where $p < q$. The
hunter taking $\varphi(x_q)$ stays there, while the other hunter moves from $\varphi(x_p)$ to the
rightmost external neighbor $\varphi(x_r)$ of $\varphi(x_q)$. By
Lemma~\ref{lem:outerstring_external}, the robber cannot move from the right of $\varphi(x_q)$ to the
left of it. Therefore, the robber either appears between $\varphi(x_q)$ and $\varphi(x_r)$, or he is
on the right of $\varphi(x_r)$. Since the sequence of external strings is finite, the hunters either
capture the robber, or he is confined by some $\varphi(x_p)$ and $\varphi(x_q)$, so the phase ends
after finitely many steps.

Since there are only finitely pairs of intersecting strings nested in each other, the strategy
proceeds in finitely many phases and the robber is captured.\qed
\end{proof}

Similarly as in Section~\ref{sec:ifa}, we can show that the strategy requires at most a linear
number of moves.

\section{Guarding Shortest Paths and Curves in String Graphs} \label{sec:guarding_path}

In this section, we build a crucial tool for designing our strategy to capture the robber using 15
cops in any string graph. The main result shows that 5 cops are able to guard a shortest curve in a
string representation together with the strings intersecting it.

\heading{Guarding Shortest Paths.}
We recall a classical lemma of Aigner and Fromme~\cite{aigner_fromme}:

\begin{lemma}[\cite{aigner_fromme}]\label{lem:guarding_shortest_path}
Let $P$ be a shortest path between a pair of vertices of a graph~$G$. Then a single cop has a strategy to guard $P$, after a finite number of
initial moves.
\end{lemma}

In~\cite{aigner_fromme}, this result is  essential to prove that the maximum cop number of planar
graphs is three.  The idea is that one can cut the planar graph by protecting several shortest
paths. Consider a planar embedding.  The strategy is to protect two shortest paths
$P_1$ and $P_2$ from $u$ to $v$ such that the robber is confined to the subgraph $D$ between $P_1$
and $P_2$. A third shortest path $P_3$ in $D$ is chosen and guarded by the third cop. The robber
has to choose one of the smaller subgraphs $D'$ of $D$ to which he is confined. It is shown that
$N(D')$ can be guarded by just two paths, so one of the cops can be freed and the strategy can be
iterated.

\heading{Guarding Retracts.}
There is the following simple generalization of Lemma~\ref{lem:guarding_shortest_path}.  It was
stated in~\cite{retracts} (in a different form) and we believe that this statement should be more
known. A \emph{retract} from $G = (V_G, E_G)$ to an induced subgraph $H = (V_H, E_H)$ is a map $f\colon
V_G \to V_H$ such that $f(v)=v$ for all $v\in V_H$, and for every $uv \in E_G$ either $f(u) f(v) \in
E_H$, or $f(u)=f(v)$.

\begin{lemma}\label{lem:guarding_retract}
Let $H$ be a retract of $G$. Then $\cn(H)$ cops have a strategy in $H$ to position one of them, in
finite number of steps.  After the positioning, the cop can guard $H$ while the remaining ones are
free.
\end{lemma}

\begin{proof}
The strategy for $\cn(H)$ cops plays on $H$ as if a robber standing on $r \in V(G)$ is placed on
$f(r)$. By the definition of the retract, $f(r)$ moves by the distance at most 1 with each turn of
the game. Therefore, $\cn(H)$ cops have a strategy to ``capture'' $f(r)$ in finitely many turns.
The cop standing on $f(r)$ can then follow the robber, to be always on $f(r)$ for the current
robber's position $r$. If the robber steps on $V(H)$, he is immediately captured by the cop. Note
that the remaining cops are no longer required.\qed
\end{proof}

This implies Lemma~\ref{lem:guarding_shortest_path} since a shortest path is a retract and paths
have the cop number equal to 1.

\heading{Guarding Neighborhoods of Shortest Paths.}
We want to apply a similar idea to string graphs. Unfortunately, guarding a shorting path $P$ is not
sufficient to prevent the robber to move from one side of $P$ to the other one.  We need a stronger
tool to geometrically restrict the robber.  We show that five cops are sufficient to guard $N[P]$
which prevents to the robber to use any string crossing the protected path; see
Fig.~\ref{fig:guarding_strings}.

Before stating the lemma, we add another definition.  Suppose that the robber is confined by the
strategy to $D \subseteq V$. We say that a path $P$ in $G$ is shortest \emph{relative to $D$}, if it
is shortest in $G[P \cup D]$. The path does not have to be shortest in $G$, it just have to be
shortest with respect to the robber's confinement to $D$:

\begin{lemma}\label{lem:guarding_shortest_path_neighbourhood}
Let $P$ be a shortest path relative to $D\subseteq V$ and let the robber be confined to $D$.
Then five cops have a strategy to guard $N[P]$, after a finite number of initial moves.
\end{lemma}

\begin{proof}
It follows from Lemma~\ref{lem:guarding_shortest_path} applied to $G[D \cup P]$ that one cop, called
the \emph{sheriff}, has a strategy to guard $P$ since the robber can only move in $D$. The four
additional cops, called the \emph{deputies}, follow the sheriff and stand at neighboring vertices of
$P$. More precisely, suppose that the path $P$ consists of the vertices $p_0,p_1,\dots,p_k$.
When the sheriff stands at $p_i$, the deputies stand at $p_{i-2}$, $p_{i-1}$, $p_{i+1}$ and
$p_{i+2}$ (with the convention that $p_{-1}$ and $p_{-2}$ here refer to the 
vertex $p_0$, and $p_{k+1}$ and $p_{k+2}$ refer to $p_k$).
As the sheriff moves along the 
path according to the strategy, the deputies follow him.
The initial setup procedure is analogous to the one in Lemma~\ref{lem:guarding_shortest_path}.

Assume that the robber moves to a vertex $r$, and suppose that $r$ is adjacent to a vertex $q$,
which is adjacent to $p_i$. Then the strategy necessarily moves the sheriff to one of the vertices
$p_{i-2},\dots,p_{i+2}$, since otherwise the robber could step on $P$ in two moves without being
immediately captured by the sheriff, contradicting the properties of the strategy. Therefore, after
the cops' move, there is at least one cop on $p_i$, and so if the robber moves to $q$, he is
captured immediately.\qed
\end{proof}

To guard $N[P]$ with the cops moving only on $P$, five cops are necessary as shown in
Fig.~\ref{fig:five_neccessary}b.  When we say that five cops \emph{start guarding a path}, we do not
explicitly mention the initial time required to position them onto the path and assume that the
strategy waits for enough turns.

Unfortunately, the result of Lemma~\ref{lem:guarding_shortest_path_neighbourhood} cannot be
straightforwardly extended to retracts. Even if the retract has bounded degree (so the number of
deputies required to guard the vertices in distance at most 2 is bounded), it is not possible to
move the deputies together with the sheriff in the required way.

\begin{figure}[t!]
\centering
\includegraphics{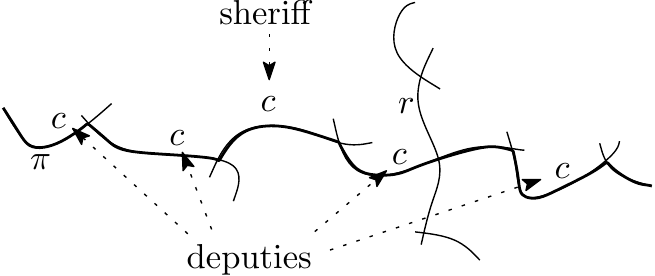}
\hskip 10mm
\includegraphics{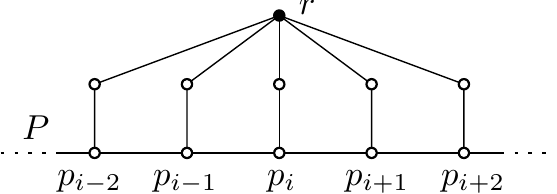}
\caption{(a) To guard a shortest curve $\pi$ defined by a path $P$, five cops guard 
consecutive strings of $\pi$. A string $r$ crossing $\pi$ may not be in $P$, but it belongs to $N[P]$. 
(b) The necessity of five cops to guard $N[P]$. With the robber standing on $r$, there
needs to be a cop on each of the vertices $p_{i-2},\dots,p_{i+2}$, otherwise the robber could
safely move to $N[P]$.}
\label{fig:guarding_strings}
\label{fig:five_neccessary}
\end{figure}

\heading{Guarding Shortest Curves.}
Our strategy for string graphs is geometric, based on string representations. To simplify its
description, we introduce the concept of shortest curves as particular curves through the string
representation of some shortest path.

Let $G$ be a string graph together with a fixed string representation $\varphi$.  Suppose that the
robber is confined to $D\subseteq V$ and let $P$ be a shortest path relative to $D$ from $u$ to $v$.
Suppose that we choose and fix two points $A \in \varphi(u)$ and $B \in \varphi(v)$.  Let $\pi
\subseteq \varphi(P)$ be a curve from $A$ to $B$ such that for every $p \in P$, the curve $\pi$ has
a connected intersection with $\varphi(p)$, and these intersections are ordered on $\pi$ in the same order
as the vertices of $P$. We call $\pi$ a \emph{shortest curve of $P$, relative to $D$} with endpoints
$A$ and $B$. A curve $\pi$ is called a \emph{shortest curve relative to $D$} if it is a shortest
curve of some shortest path relative to $D$. We may omit $D$ when it is clear from the context. 

The shortest path corresponding to a shortest curve $\pi$ is uniquely defined by the sequence of
strings whose intersection with $\pi$ has non-zero length.  By \emph{guarding a shortest curve
$\pi$}, we mean guarding $N[P]$ of the corresponding shortest path $P$. The \emph{length} of $\pi$
is the number of its strings; we note that its Euclidean length plays no role.

\begin{corollary}\label{lem:guarding_shortest_curve}
Let the robber be confined to $D$ and let $\pi$ be a shortest curve relative to $D$  in a string
representation.  Then five cops can prevent the robber from entering any string intersecting $\pi$,
after a finite number of initial moves.
\end{corollary}

\begin{proof}
Let $P$ be the shortest path defining $\pi$. By guarding $N[P]$, five cops prevent the robber from
entering strings intersecting $\pi$. See Fig.\ref{fig:guarding_strings}a for illustration.\qed
\end{proof}

\begin{observation} \label{lem:subcurve}
Any sub-curve of a shortest curve relative to $D$ is a shortest curve relative to $D$.\qed
\end{observation}

\section{Capturing Robber in String Graphs}\label{sec:string}

In this section, we show that the maximum cop number of string graphs is at most 15.  Our strategy
is inspired by the strategy for 3 cops in planar graphs~\cite{aigner_fromme}.  The key difference is
that we use Lemma~\ref{lem:guarding_shortest_path_neighbourhood} instead of
Lemma~\ref{lem:guarding_shortest_path}, so we require 5 cops for each shortest path instead of 1.
Therefore, our strategy requires $3 \cdot 5 = 15$ cops.

\heading{Segments, faces and regions.}
Consider a set $\mathcal C$ of curves/strings in $\R^2$.
The topological arc-connected components of $\R^2 \setminus \mathcal C$ are called \emph{faces} and
their topological closures are \emph{closed faces}; every face is an open set.
We assume that the number of intersections of $\calC$ is finite, so the number
of faces is also finite.

A \emph{segment} of a curve $\pi \in \mathcal C$ is a maximal arc-connected subset of $\pi$ not containing any intersection with
another curve in $\calC$.
The number of segments is also finite. A \emph{region} is a closed subset of $\R^2$
obtained as a closure of a union of some of the faces.

Consider a string representation $\varphi$.  For $X \subseteq \R^2$, we denote the topological
closure of $X$ by $\closure X$, the topological interior by $\interior(X)$, and the boundary $\bo X
= \closure X \setminus \interior(X)$. We say that a vertex $v$ is \emph{contained in $X$} if $\varphi(v)
\subseteq \interior(X)$.  We denote the subgraph of $G$ induced by the vertices contained in $X$ by
$G_X$.  Two curves sharing only their endpoints are said to to be \emph{internally disjoint}. 

\begin{lemma} \label{lem:curves_confinment}
Let $\pi_1$ and $\pi_2$ be two internally disjoint curves with endpoints $a$ and $b$, let $F$ be the
closed face of $\R^2\setminus(\pi_1\cup\pi_2)$ containing the string with the robber, and let $D$ be
the connected component of $G_F$ containing the robber.  If $\pi_1$ and $\pi_2$ are shortest curves
with respect to $D$, each guarded by five cops, then the robber is confined to $D$.
\end{lemma}

\begin{proof}
It follows from Corollary~\ref{lem:guarding_shortest_curve} applied to $\pi_1$ and $\pi_2$.\qed
\end{proof}

Additionally, below we use the following topological lemma.

\begin{lemma}\label{lem:zigzag}
Let $\pi_1$ and $\pi_2$ be two internally disjoint simple curves from $a$ to $b$, where $a \ne b$.  Let $F$ be
a closed face of $\R^2 \setminus (\pi_1\cup\pi_2)$. Let $\pi_3 \subseteq F$ be a simple curve from
$a$ to $b$, going through at least one inner point of~$F$. Then every face $R$ of $F \setminus (\pi_1 \cup
\pi_2 \cup \pi_3)$ is bounded by two simple internally disjoint curves $\pi_i'$ and $\pi_3'$, where
$\pi_i'\subseteq\pi_i$ for some $i\in\{1,2\}$ and $\pi_3'\subseteq\pi_3$.
\end{lemma}

\begin{proof}
Without loss of generality, we may assume that $F$ is the inner face of $\R^2 \setminus (\pi_1 \cup \pi_2)$,
otherwise we can apply the circular inversion.  We know that $R$ is an open arc-connected set by
definition, so $\bo R$ is a simple closed Jordan curve.

We first establish that $\bo R \subseteq (\pi_i \cup \pi_3)$ for some $i \in \{1,2\}$.  Observe
that $\bo R \subseteq \pi_j$ would imply that $\pi_j$ is not a simple curve.  There is a point $r_3
\in (\bo R \cap \pi_3) \setminus (\pi_1 \cup \pi_2)$. The reason is that otherwise we
would have $\bo R = \pi_1 \cup \pi_2$, so $\closure R = F$, which contradicts that $\pi_3$
intersects $\interior(F)$.

We argue that it is not possible that there exist both $r_1 \in (\bo R \cap \pi_1) \setminus (\pi_2
\cup \pi_3)$ and $r_2 \in (\bo R \cap \pi_2) \setminus (\pi_1 \cup \pi_3)$. If both would exists,
there would be a curve $\pi' \subseteq (R \cup \{r_1,r_2\})$ from $r_1$ to $r_2$ separating $a$ from
$b$ in $F$, not intersecting $\pi_3$. By Jordan curve theorem, this contradicts that $\pi_3$
is a curve from $a$ to $b$ through $F$; see Fig.~\ref{fig:zigzag}a. However, we
necessarily have one such $r_i$, for $i \in \{1,2\}$. Without loss of generality, we assume that
$r_1$ exists and no $r_2$ exists, so $\bo R \subseteq \pi_1 \cup \pi_3$.

\begin{figure}[t!]
\centering
\includegraphics{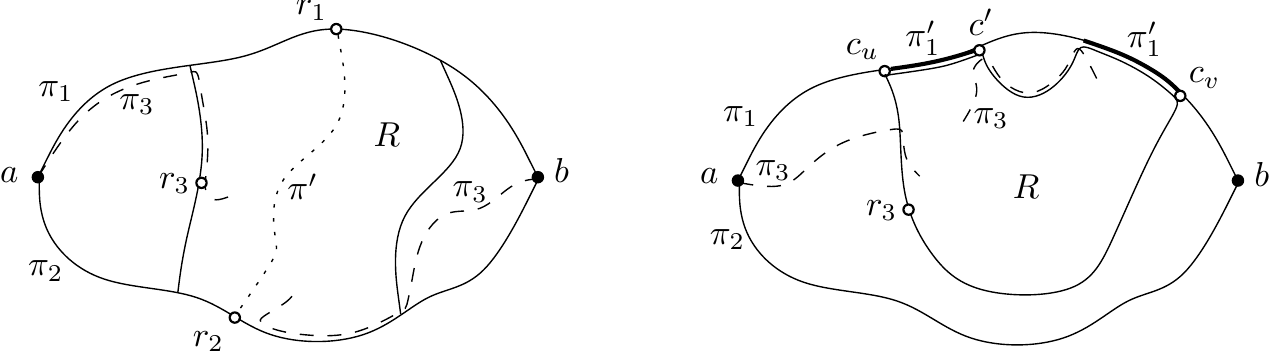}
\caption{(a) The proof of Lemma~\ref{lem:zigzag}: there cannot be both $r_1$ and $r_2$ as in the
proof. (b) A disconnected $\pi_1'$ (bold line) implies that $\pi_3$would cross $R$.}
\label{fig:zigzag}
\end{figure}

Let $\pi_1'= \closure{(\bo R \cap \pi_1) \setminus \pi_3}$. Let $c_u$ be the first point of $\pi_1'$
going along $\pi_1$ from $a$ to $b$ and $c_v$ last such point. If $\pi_1'$ was not
connected, let $c' \ne c_u,c_v$ be an endpoint of one segment of $\pi_1'$; 
necessarily, $c' \in \pi_3$ since $\bo R \subseteq \pi_1 \cup \pi_3$.  However, it is not possible
since $\pi_3$ would have to contain a point of $R \cup (\R^2 \setminus F) \cup \pi_1'$, as shown in
Fig.~\ref{fig:zigzag}b, contradicting the definitions of $\pi_3$ and $R$. Therefore $\pi_1'$ is
connected and we can take $\pi_3' = \closure{\bo R \setminus \pi_1'}$, getting a connected curve
$\pi'_3 \subseteq \pi_3$.\qed
\end{proof}

\heading{Restricted Graphs and Strategies.}
Given a closed region $R \subseteq \R^2$, let \emph{$G$ restricted to $R$}, denoted $G|_R$, be the
intersection graph of the strings of $\varphi \cap R$.  This restriction may remove vertices
(represented by the strings outside $R$), may remove edges (intersections outside $R$) and may split
vertices whose strings leave and reenter $R$ at least once; every arc-connected part of $\varphi(v)
\cap R$ forms a new vertex $v_i$. The newly obtained vertices $v_i$ are called the \emph{splits} of
$v$.  The graph $G|_R = (V|_R, E|_R)$ is again a string graph with its representation, denoted by $\varphi|_R$, directly
derived from $\varphi$.  Note that this operation preserves the faces and strings in $\interior(R)$
and all representation properties assumed above, namely the vertex set of $G|_R$ is finite.  Also,
the number of segments does not increase.

\begin{lemma} \label{lem:restricted_shortest_curve}
Let $R$ be a region. If $\pi$ is a shortest curve relative to $D$, and $\pi' \subseteq R$ is a
sub-curve of $\pi$, then $\pi'$ is a shortest curve relative to $D \cap V|_R$ in
$G|_R$ and $\varphi|_R$.
\end{lemma}

\begin{proof}
Observe that the underlying path $P'$ of $\pi'$ is preserved, and if any $p \in P'$ is split in
$G|_R$, we use $p_i$ intersecting $\pi'$. The rest follows from Observation~\ref{lem:subcurve}, and
the fact that no path is shortened in $G|_R$.\qed
\end{proof}

We now show that certain strategies for a restricted graph can be used in the original graph.

\begin{lemma} \label{lem:restricted_strategy}
Let $R$ be a region such that the robber is confined to $G_R$.  Suppose that there
exists a strategy $\mathcal S'$ capturing the robber in $G|_R$, confining him for the entire
strategy to $G_R$.  Then there exists a strategy $\mathcal S$ for the same number of cops capturing
the robber on $G$, if the robber is initially confined to $G_R$.
\end{lemma}

\begin{proof}
The strategy $\mathcal S$ proceeds as $\mathcal S'$, with the following exception.
When $\mathcal S'$ moves a cop to a split $v_i\in V_{G|_R}$ of $v\in V_G$, the strategy $\mathcal S$
move this cop to $v$; note that this move is always possible. It is key that the robbers choices
in $G_R$ are not extended, so he is confined to it by $\mathcal S$ and captured.
\qed
\end{proof}

\heading{Proof of Theorem~\ref{thm:bounded_max_cn}(iii).}
We are ready to prove that the maximum cop number of string graphs is at least 3 and at most~15.

\begin{proof}[Theorem~\ref{thm:bounded_max_cn}({\rm iii})]
The lower bound of 3 cops follows from the graph in Fig.~\ref{fig:mrizka}. It remains to argue that
there exists a strategy using 15 cops. Our strategy proceeds in phases, monotonously shrinking the
confinement of the robber. In the beginning of each phase, the robber is confined to $D\subseteq V$
either (A) by a single cop guarding a cut-vertex separating $D$ from the rest of the graph, or (B)
by ten cops guarding two shortest curves forming a simple (non-self-intersecting) cycle surrounding
$D$. In each phase, we decrease the number of vertices in $V$ or $D$, so the robber is caught after
finitely phases.

Let $B$ be the union of the currently guarded paths 
and vertices; by Lemma~\ref{lem:guarding_shortest_path_neighbourhood}, if the robber moves to $N[B]$, he
is captured. Let $D$ be the component of $G \setminus N[B]$ containing the vertex with the robber,
and let $Q=N[B] \cap N[D]$. Since our strategy confines the robber to $D$ for the rest of the game,
we can leave out the remaining vertices and assume that $V=D \cup Q \cup B$. Let $s$ be the number
of segments of $\varphi$.

\begin{claim} \label{claim:string_proof}
Let $V = D \cup Q \cup B$, the robber stands on $r \in D$, and one of the following holds:
\begin{packed_itemize}
\item[(A)] $|B|=1$ and 1 cop guards a vertex $c \in B$.
\item[(B)] $|B| \ge 2$, 10 cops guard two shortest curves $\pi_1$ and $\pi_2$ relative to $D$ between points
$a$ to $b$ such that $\pi_1 \cup \pi_2$ forms a simple cycle, and additionally
$G=G|_F$ 
where $F$ is the closed face of $\R^2\setminus(\pi_1\cup\pi_2)$ containing $\varphi(r)$.
\end{packed_itemize}
Then 15 cops have a strategy to capture the robber.
\end{claim}

\begin{proof}[Claim]
We prove this claim by induction on $s$ and $|D|$. The claim obviously true when $s \le 1$ and $|D|
= 0$. The strategy proceeds differently according to which of (A) and (B) is satisfied.

\emph{Case (A).} If $Q=\{q\}$, then move the cop guarding $c$ to start guarding $q$.  Let $G'=G -
c$, we further leave out the irrelevant vertices, so $V'=D'\cup Q'\cup \{q\}$ as above. The rest follows
from the induction hypothesis, with the assumption (A), applied to $G'$ with $s' \le s$ and
$D'\subsetneq D$.

Let $Q = \{q_1,\dots,q_k\}$ for $k \ge 2$. Let $a_i$ be a point of $\varphi(c) \cap \varphi(q_i)$.
We choose $\pi_1$ be a shortest curve in $\varphi(V \setminus \{c\})$ between some $a_i$ and $a_j$,
and $\pi_2$ be the subcurve of $\varphi(c)$ between $a_i$ and $a_j$.  Without loss of generality,
$\pi_1 \cup \pi_2$ forms a simple cycle; if not, we can shorten it by choosing different points $a_i$ and
$a_j$.

We start guarding $\pi_1$ and $\pi_2$ with 10 cops.  Let $F'$ be the closed face of
$\R^2\setminus(\pi_1\cup\pi_2)$ containing a string with the robber.  We denote $G' = G |_{F'}$, we
leave out the irrelevant vertices, so $V'=D'\cup Q'\cup B'$ as above.  We use the induction
hypothesis, with the assumption (B), applied to $G'$ for $s' \le s$ and $D' \subsetneq D$.
By Lemma~\ref{lem:restricted_strategy}, the strategy on $G'$ from the induction hypothesis implies a
strategy on $G$.

\emph{Case (B).} First suppose that there exists no shortest curve in $\varphi$ between $a$ and $b$
intersecting $\interior(F)$. By Menger theorem, there must be a cut-vertex $c \in B \cup Q$
separating $D$ from $B$. Our strategy guards $c$ with one cop, and then stops guarding $B$. Let $G'=
G\bigl[(V \setminus B) \cup \{c\}\bigr]$, leaving out the irrelevant vertices, so $V'=D' \cup Q'\cup
\{c\}$ as above.  The rest follows from the induction hypothesis, with the assumption (A), applied
to $G'$ with $s' < s$ and $D'\subseteq D$.

Otherwise, let $\pi_3$ be a shortest curve relative to $D$ in $\varphi$ from $a$ to $b$ intersecting
$\interior(F)$. The strategy starts guarding $\pi_3$ with the five free cops. Then, let $F'$ be the
closed face of $\R^2\setminus(\pi_1\cup\pi_2\cup\pi_3)$ containing the string on which the robber
stands. By Lemma~\ref{lem:zigzag}, we have that $\bo F' = \pi'_i \cup \pi'_j$ where $\pi'_i$ is a
subcurve of $\pi_i$, $\pi'_j$ is a subcurve of $\pi_j$, and $\pi_i'\cup\pi_j'$ form a simple cycle.
We free the five cops stop guarding $\pi_k$, where $k \ne i,j$, and we restrict the guarding of
$\pi_i$ and $\pi_j$ to $\pi_i'$ and $\pi_j'$, which are shortest curves by
Observation~\ref{lem:subcurve}.

Let $G'=G|_{F'}$, leave out the irrelevant vertices, so $V'=D'\cup Q'\cup B'$ as above. The rest
follows from the induction hypothesis, with the assumption (B), applied to $G'$ with
$s' \le s$ and $D'\subsetneq D$. By Lemma~\ref{lem:restricted_strategy}, the strategy on $G'$
from the induction hypothesis implies a strategy on $G$.\claimqed
\end{proof}

The theorem follows by guarding an arbitrary vertex $c$ with one cop, so $B=\{c\}$. We leave out the
irrelevant vertices, so $V' = D \cup Q \cup B$. We use Claim~\ref{claim:string_proof} with the
assumption (A) for $G'=G|_{V'}$.\qed
\end{proof}

\section{Capturing Robber in String Graphs on Bounded Genus Surfaces}\label{sec:bounded_genus}

In this section, we generalize the results of the previous section to graphs having a string
representation on a fixed surface, and we prove Theorem~\ref{thm:bounded_max_cn}(iv).

\heading{Definitions.}
We assume familiarity with basic topological concepts related to curves on surfaces, such as genus, Euler genus,
non-contractible closed curves, the fundamental group of surfaces and graph embedding properties.  A
suitable treatment of these notions can be found in \cite{prasolovelements,MoTh}. 

A \emph{walk} $W$ in a graph $G$ is a sequence $w_0, w_1,\dots w_k$ of vertices where $w_i$ and
$w_{i+1}$ adjacent; repetitions of vertices and edges are allowed.  A walk is called \emph{closed}
if $w_0=w_k$. Let $|W|$ denote the length $k$ of $W$.  For walks $W=w_0,w_1,\dotsc,w_k$ and
$W'=w'_0,w'_1,\dotsc,w'_\ell$ with $w_k=w'_0$, we denote the concatenation by $W+W' =
w_0,w_1,\dotsc,w_k,w'_1,w'_2,\dotsc,w'_\ell$. Let $-W$ be the reversal of $W$ and let
$W_1-W_2=W_1+(-W_2)$.

A curve $\pi$ is a continuous function from the interval $[0,1]$ to the surface, and it is
\emph{closed} if $\pi(0) = \pi(1)$. The concatenation of curves $\pi_1+\pi_2$ is defined naturally
whenever $\pi_1(1)=\pi_2(0)$, and similarly $-\pi$ is the reversal and $\pi_1-\pi_2=\pi_1+(-\pi_2)$.
We use the following topological lemma, following from the properties of the fundamental group;
see~\cite{prasolovelements}.

\begin{lemma}[\cite{prasolovelements}] \label{lem:three-curves}
Let $\pi_1$, $\pi_2$ and $\pi_3$ be three curves on a surface $\mathbf{S}$ from $a$ to $b$.
If the closed curve $\pi_1-\pi_2$ is non-contractible, then at least one of 
$\pi_1-\pi_3$ and $\pi_2-\pi_3$ is non-contractible.
\end{lemma}

Consider a string representation $\varphi$ of $G$ on a surface~$\mathbf{S}$.  We represent the
combinatorial structure of $\varphi$ by an auxiliary multigraph $A(\varphi)$ embedded on
$\mathbf{S}$ defined as follows.  The vertices of $A(\varphi)$ are the endpoints of the strings of
$\varphi$ and the intersection points of pairs of strings of~$\varphi$. The edges of $A(\varphi)$
correspond to segments of strings of $\varphi$, i.e., to subcurves connecting 
pairs of vertices appearing consecutively
on a string of~$\varphi$. By representing $\varphi$ by $A(\varphi)$, we can use the well-developed
theory of graph embeddings on surfaces. 

\heading{Walks Imitating Non-contractible Curves.}
We introduce a relation between a walk in $G$ and a curve on $\mathbf{S}$, allowing us to easily
transition between the two.  We say that a walk $W=w_0,w_1,\dots w_k$ in $G$ \emph{imitates} a curve
$\pi\subseteq \varphi(G)$ on the surface $\mathbf{S}$ if $\pi$ can be partitioned into a sequence of
consecutive subcurves $\pi_0,\pi_1,\dotsc,\pi_k$ of positive length such that
$\pi=\sum_{i=0}^{k}\pi_i$ and $\pi_i\subseteq\varphi(w_i)$. A closed walk $W$ \emph{imitates a
non-contractible curve} if there is a non-contractible curve $\pi\subseteq\varphi(G)$ imitated
by~$W$.

\begin{lemma}\label{lem:cut-surface}
Let $\varphi$ be a string representation of $G$ on an orientable (resp. non-orientable) surface $\mathbf{S}$ of genus $g>0$ (resp. Euler genus $g'>0$) and let $W$
be a closed walk in $G$ imitating a non-contractible curve. Then every connected component of $G \setminus N[W]$ has a string representation on a surface of genus at most $g-1$ (resp. Euler genus at most $g'-1$).
\end{lemma}

\begin{proof}
Note that the proof and te arguments are the same for orientable genus and Euler genus.

If $A(\varphi)$ has an embedding on a surface of genus $g-1$ (resp. Euper genus $g'-1$), then $G$ has a string representation
on this surface and we are done. Suppose then that this is not the case, i.e., $A(\varphi)$ is a graph of genus $g$ (resp Euler genus $g'$). Therefore its embedding on $\mathbf{S}$ is a 2-cell embedding, i.e., every face of $\mathbf{S}-\varphi$ is homeomorphic to a disk.

Let $\pi$ be the non-contractible curve imitated by~$W$.  The curve $\pi$ traces a closed walk $W'$
in~$A(\varphi)$. Since $\pi$ is non-contractible, $W'$ contains a non-contractible simple cycle $C$
of~$A(\varphi)$.  By standard results on 2-cell embeddings (see \cite[Chapter 4.2]{MoTh}), the genus (resp. the Euler genus)
of every connected component of $A(\varphi) \setminus C$ is strictly smaller than the genus of
$A(\varphi)$.

Let $\varphi'$ be the string representation $\varphi|_{V \setminus N[W]}$.  The auxiliary multigraph
$A(\varphi')$ is a subgraph of $A(\varphi) \setminus C$, and hence each of its connected components
has an embedding on a surface of genus $g-1$ (resp. Euper genus $g'-1$). This embedding corresponds to a string representation
of a connected component of $G \setminus N[W]$ on a surface of genus~$g-1$ (resp. Euper genus $g'-1$).\qed
\end{proof}

\begin{lemma}\label{lem:non-contractible-imitated}
If $G$ has no string representation in the plane, then for every string representation $\varphi$ of
$G$ on a surface $\mathbf{S}$ there is a closed walk $W$ in $G$ imitating a non-contractible curve.
\end{lemma}

\begin{proof}
Since $A(\varphi)$ is not planar, the embedding of $A(\varphi)$ contains a non-contractible cycle
(see~\cite[Chapter 4.2]{MoTh}), which corresponds to a non-contractible curve on~$\mathbf{S}$. This
curve is imitated by a closed walk $W$ of~$G$.\qed
\end{proof}

\begin{lemma}\label{lem:three-walks}
Let $\varphi$ be a string representation of $G$ on a surface~$\mathbf{S}$, 
let $u,v\in V$ be two vertices, and let $W_1$, $W_2$, $W_3$ be three walks from $u$ to $v$.
If $W_1-W_2$ imitates a non-contractible closed curve, then at least one of $W_1-W_3$ and
$W_2-W_3$ imitates a non-contractible closed curve.
\end{lemma}

\begin{proof}
Let $\pi_{12}$ be a non-contractible closed curve imitated by $W_1-W_2$. Looking at the 
consecutive subcurves of $\pi_{12}$ corresponding to vertices of $W_1-W_2$, we 
have that
$\pi_{12}=\pi_1-\pi_2$ with $\pi_1$ imitated by $W_1$ and $\pi_2$ imitated by $W_2$. Let
$x\in\varphi(u)$ be the first point of $\pi_1$ and $y\in\varphi(v)$ be the last point of $\pi_1$.

Now let $\pi_3$ be any $x$-$y$ curve imitated by $W_3$ and observe 
that $\pi_1-\pi_3$ is imitated by $W_1-W_3$ and $\pi_2-\pi_3$ is imitated by $W_2-W_3$.
By Lemma~\ref{lem:three-curves}, at least one of $\pi_1-\pi_3$ and 
$\pi_2-\pi_3$ is non-contractible
and the lemma follows.\qed
\end{proof}

\begin{lemma}\label{lem:guard-shortest-loop}
On a graph $G$ with a string representation $\varphi$ on a surface $\mathbf{S}$ and a shortest closed walk
$W$ imitating a non-contractible curve, 10 cops have a strategy
to guard $N[W]$ after a finite number of initial moves.
\end{lemma}

\begin{proof}
If $|W|\leq 10$, the cops may occupy every vertex of $W$ for the rest of the game and we are done.
Otherwise we divide $W$ into two almost-equally long walks $W_1$, $W_2$, where $W_i$ is from $u_i$
to $v_i$, where $u_1u_2$, $v_1v_2$ are edges. We have $|W_1|\geq|W_2|\geq|W_1|-1$, such that
$W=W_1+v_1v_2-W_2-u_1u_2$; see Fig.~\ref{fig:closed-walk-split}. Note that $|W|=|W_1|+|W_2|+2$.

\begin{figure}[t!]
\centering
\includegraphics{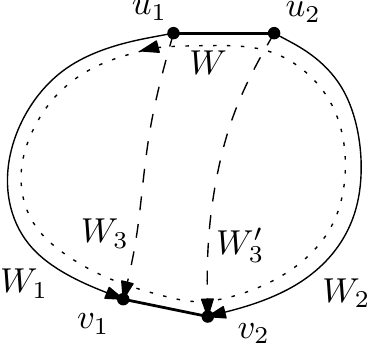}
\caption{The situation in the proof of Lemma~\ref{lem:guard-shortest-loop}.}
\label{fig:closed-walk-split}
\end{figure}

We claim that both $W_1$ and $W_2$ are shortest paths in $G$. If $W_1$ is not a 
shortest path, let $W_3$ be a shortest path from $u_1$ to $v_1$, so $|W_3|<|W_1|$.
Then both closed walks $W_1-W_3$ and $W_3+v_1v_2-W_2-u_1u_2$ would be shorter than $W$:
$$|W_1-W_3|=|W_1|+|W_3|<2|W_1|\leq|W_1|+|W_2|+1<|W|,$$
$$|W_3+v_1v_2-W_2-u_1u_2|=|W_2|+|W_3|+2<|W_1|+|W_2|+2=|W|.$$
By Lemma~\ref{lem:three-walks}, at least one of them is non-contractible which contradicts the
assumption. Similarly, there exists no path $W'_3$ from $u_2$ to $v_2$ with $|W'_3|<|W_2|$.

Therefore we may use Lemma~\ref{lem:guarding_shortest_path_neighbourhood} (with $D=V$) to guard
$N[W_1]$ and $N[W_2]$ with ten cops.\qed
\end{proof}

\heading{Proof of Theorem~\ref{thm:bounded_max_cn}(iv,v).}
We are ready to prove that the maximum cop-number of \genusstr{g} graphs
is at least  ${1 \over 2}g^{1 \over 3}$ and at most~$10g+15$, and of \eulergenusstr{g'} graphs is at
least ${1 \over 2}{g'}^{1 \over 3}$ and at most~$10g'+15$.

\begin{proof}[Theorem~\ref{thm:bounded_max_cn}({\rm iv,v})]
For the lower bound, consider the incidence graph $G$ of a projective plane of order
$k=\lceil{1 \over 2}g^{1 \over 3}\rceil$. It has at most than $2k^3$ edges, so its genus is at most
$g$. It is $(k+1)$-regular graph with girth 6, so $\cn(G) \ge k+1 \ge {1 \over 
2}g^{1 \over 3}$ (in fact, $\cn(G)=k+1$, as shown by Pra\l at~\cite{pralat}). A similar lower bound
works for non-orientable surfaces since $g' \le 2g$.

\begin{figure}[b!]
\centering
\includegraphics{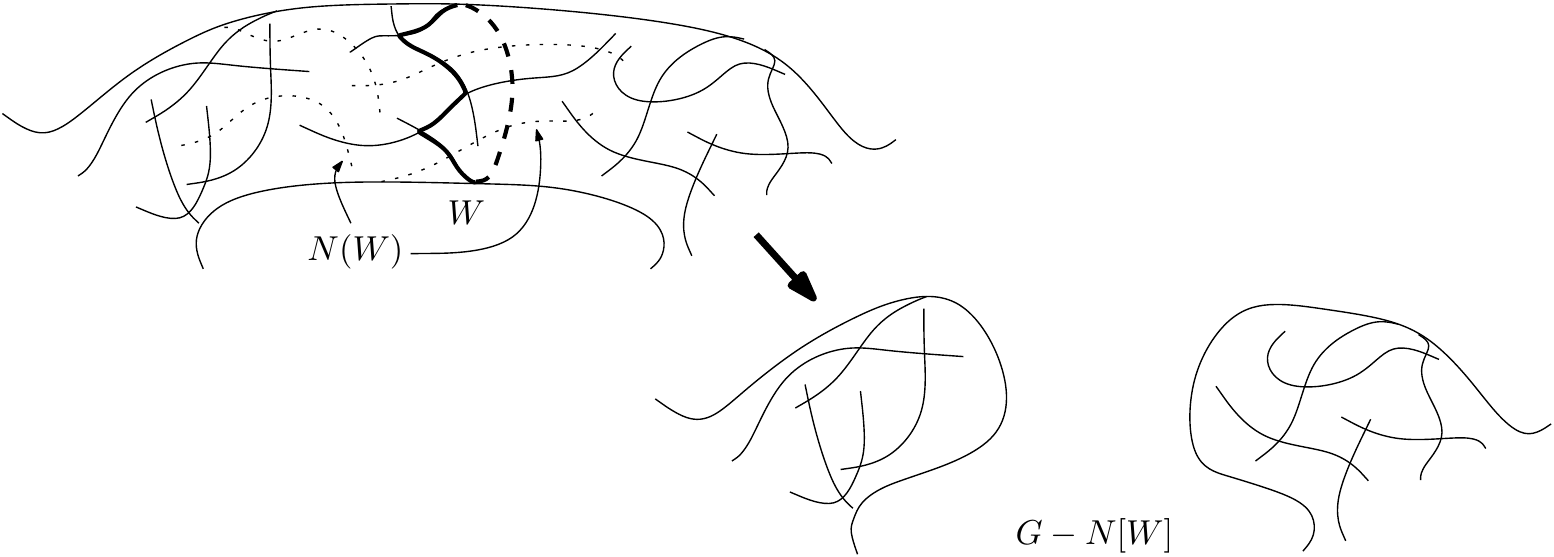}
\caption{Cutting a surface handle after guarding $N[W]$ where $W$ is a closed walk imitating a
non-contractible curve. Note that our proof is general and works for any non-contractible curve.}
\label{fig:handle-split}
\end{figure}
For the upper bound, we proceed by induction on the genus $g$ in case of orientable surfaces, resp. by the Euler genus $g'$ in case of non-orientable surfaces. The proof here is the same for both parameters. Note that a non-orientable surface may become orientable after removing a non-contractible curve (as in Lemma~\ref{lem:cut-surface}), the proof for Euler genus does not depend on the underlying surface being non-orientable.

The cases $g=0$ and $g'=0$ are proved by Theorem~\ref{thm:bounded_max_cn}(iii). Suppose now that $g>0$ (resp. $g'>0$) and fix a string representation $\varphi$ of $G$ on a surface of genus~$g$ (resp. Euler genus $g'$). Let $W$ be a shortest closed walk in $G$ imitating a
non-contractible curve. By Lemma~\ref{lem:guard-shortest-loop}, 10 cops prevent the robber from entering~$N[W]$ till the end of the game; see Fig.~\ref{fig:handle-split}.

Thus, after a finite number of moves the robber will remain confined to a connected component $G'$ of $G \setminus N[W]$. By Lemma~\ref{lem:cut-surface}, $G'$ has a string representation on a surface of
genus at most $g-1$ (resp. Euler genus $g'$).  By the induction hypothesis, $15+10(g-1)$ cops (resp. $15+10(g'-1)$ cops) have a strategy to capture the robber on $G'$ and, together with 10 cops gurading $W$, also on $G$.\qed
\end{proof}

\section{Unbounded Cop Number of Intersection Graphs of Disconnected or 3-Dimensional Sets}
\label{sec:unbounded_max_cn}

In this section, we prove Theorem~\ref{thm:unbounded_max_cn} stating that the maximum cop number is
$+\infty$ even for very simple intersection classes of disconnected or 3-dimensional sets.

\heading{Cop Number of Subdivisions.}
For a graph $G=(V,E)$ and an integer $d\ge 1$, let $\Gsub{d}$ denote the graph obtained from $G$ by
replacing each edge $e=xy \in E$ by a path $P_e$ of length $d$ connecting $x$ and~$y$. In other
words, $\Gsub d$ is obtained from $G$ by subdividing each edge of $G$ by $d-1$ new vertices. The
vertices of $\Gsub{d}$ that subdivide an edge of $G$ are the \emph{subdividing vertices}, while the
vertices of $\Gsub{d}$ belonging to $G$ are the \emph{branching vertices}. The path $P_e$ is
the \emph{edge-path corresponding to~$e$}.

\begin{lemma} \label{lem:cop_number_subdivision}
For a connected graph $G=(V,E)$ and an integer $d \ge 1$, we have
$$\cn(G) \le \cn(\Gsub{d}) \le \cn(G)+1.$$
\end{lemma}

\begin{proof}
\emph{Part 1: the inequality $\cn(G) \le \cn(\Gsub{d})$.}
On $\Gsub{d}$, we consider the $(d,d)$-game which is a modification of the standard game of cops and
robber. In it, both the robber and the cops are allowed to make $d$ consecutive moves in each turn,
instead of just one move. In other words, in the $(d,d)$-game the robber and the cops are allowed to
move, in each turn, to any vertex at distance at most $d$ from their current position.
Note that if the cops have a winning strategy for the standard game on a graph $H$, then
they also have a winning strategy for the $(d,d)$-game on $H$: a $d$-fold move of the robber can be
interpreted as a sequence of $d$ simple moves,
and each can be reacted according to the winning strategy for the standard game.

We say that a vertex $\bar v$ of the graph $G$ \emph{approximates} a vertex $v$ of $\Gsub d$ if the
following holds: either $v$ is a branching vertex and $\bar v=v$, or $v$ is a subdividing vertex
belonging to an edge-path $P_e$ and $\bar v \in e$.  Notice that if $u$ and $v$ are two vertices at
distance at most $d$ in $\Gsub d$, and if $\bar u$ is a vertex of $G$ approximating $u$, then there
is a vertex $\bar v \in V$ approximating $v$ such that $\bar u$ and $\bar v$ have distance at most
one in~$G$.

Let $\mathcal S$ be a winning strategy for $k$ cops in the $(d,d)$-game on $\Gsub d$. We now
describe a winning strategy for $k$ cops playing the standard game on~$G$. Each cop playing on $G$
is identified with one cop in $\mathcal S$. When the strategy $\mathcal S$ moves the cop to a vertex
$v$ in $\Gsub d$, the corresponding cop moves to a vertex $\bar v$ of $G$ that approximates~$v$. As
argued above, this is always possible, so each move in $\Gsub{d}$ can be performed by an
approximating move in $G$. When the robber moves from $u$ to $v$ in $G$, this move can be translated
into a $d$-fold move from $u$ to $v$ in $\Gsub d$. Therefore, the strategy $\mathcal S$ can be used
to find the corresponding response on $G$, approximating the response on $\Gsub d$. Since $\mathcal
S$ captures the robber, the cops win the game in~$G$.

\emph{Part 2: the inequality $\cn(\Gsub d) \le \cn(G)+1$.} 
Suppose that $k$ cops have a winning strategy $\mathcal S$ for the standard cops and robber game on~$G$.
We use $k+1$ cops on $\Gsub d$. The first $k$ cops, called the \emph{regular cops}, are identified with
the $k$ cops of $\mathcal S$, while the remaining cop, called the \emph{tracker}, follows a special
strategy.  In the beginning of the game, the tracker follows a shortest path towards the vertex
initially occupied by the robber.  As soon as the tracker reaches a vertex previously occupied by
the robber, he only moves along the edges previously used by the robber. More precisely, if the
tracker stands in a vertex $x$, he moves through the edge that was used by the robber during his
most recent departure from~$x$.

Suppose that the robber moves from $u$ to $v$. We say that this move is a \emph{hesitant move} if
either the robber \emph{stays in a vertex} (so $u=v$), or the robber \emph{retraces an edge} (his
immediately preceding move was from $v$ to $u$). When the robber makes a hesitant move, his distance
to the tracker decreases. Thus, the tracker ensures that the robber can only make a limited
number of hesitant moves without getting captured.

The strategy of the regular cops works as follows. Recall that each regular cop corresponds to a cop in~$\mathcal S$. For simplicity, we first assume that the robber never makes a hesitant move. In the
beginning of the game, the regular cops occupy the initial positions prescribed by~$\mathcal S$.
They wait at these positions until the robber first reaches a branching vertex~$u$; the robber makes
no hesitant moves so he must reach a branching vertex within the first $d$ moves. Suppose that the
robber moves from the branching vertex $u$ to a subdividing vertex of an edge-path $P_e$, where
$e=uv$ is an edge of~$G$. Since the robber makes no hesitant moves, he moves along $P_e$ all the way
to~$v$. After the robber enters $P_e$, each cop looks up in $\mathcal S$ the prescribed response to
the robber's move from $u$ to~$v$. If $\mathcal S$ says a cop should move from a vertex $x$ to a
vertex $y$ along an edge $e'$, the corresponding cop spends $d$ moves moving from $x$ to $y$
along~$P_{e'}$. Thus, after $d$ rounds of the play, both the regular cops and the robber will again
occupy branching vertices, and the sequence of $d$ rounds corresponded to a single round
of~$\mathcal S$. The regular cops then repeat the same process, imitating the moves of~$\mathcal S$,
until the robber is caught.

It remains to deal with the robber's hesitant moves. If the robber stays in a
vertex, all regular cops stay in their vertices as well. If the robber retraces an edge, all regular cops also
retrace their last used edges. So if the robber starts moving from $u$ to $v$ along $P_e$ in
$\Gsub{d}$ and then starts moving back, the regular cops mimic him: if a regular cop moves from $x$ to $y$
along $P_{e'}$, he keeps the same distance on $P_{e'}$ to $x$ as the robber on $P_e$ to $u$. (Notice
that each change in the direction of the robber's moves on $P_e$ is a hesitant move.) Since the
robber has a limited number of hesitant moves to avoid getting captured by the tracker, the strategy
$\mathcal S$ applies and he is captured by one of the regular cops.\qed
\end{proof}

\heading{Proof of Theorem~\ref{thm:unbounded_max_cn}.}
We are ready to prove that the maximum cop number of intersection graphs of disconnected or higher
dimensional sets is $+\infty$.

\begin{proof}[Theorem~\ref{thm:unbounded_max_cn}]
As mentioned in the introduction, the class $\linegr$ of line graphs has unbounded 
cop number by results of Dudek et al.~\cite{edges}. Moreover, each line graph 
can be represented as the intersection graph of two-element subsets of the real 
line, and therefore the classes \kunitint{2} and \kint{2} contain $\linegr$ as a subclass.
It follows that these classes have unbounded cop number as well.

%

Let us now deal with geometric intersection classes of higher-dimensional objects.
The class \grid{3} (and thus also \boxicity{3}) contains $\Gsub{3}$ for all graphs $G$: the vertices
are represented by long parallel segments, say, in the direction of the $z$-axis, having pairwise
different $x$ and $y$ coordinates. Each edge is represented by an L-shape (consisting of two segments),
connecting the parallel segments representing the corresponding vertices. 
We may even assume that each segment of the representation has unit length.

For \cube{3} and \ball{3}, draw any graph $G$ in space without crossing of edges in such a way that
all edge-curves have the same length. Also ensure that around every vertex there is a ball
containing only the initial parts of the incident edge-curves and that these parts are
straight segments. Let $a$ be the minimum of all the diameters of these balls and the distances between
edge-curves outside of the balls, note that $a>0$.

Now replace every branching vertex by a cube/ball of size sufficiently
smaller than $a$ and notice that then there are disjoint tubular corridors of a positive diameter
around every edge-curve outside the vertex cubes/balls. Therefore there exists a suitable value of $d$
(depending on $G$ and the curve representation) such that $\Gsub{d}$ can be represented by chains
of sufficiently small cubes/balls within these corridors.

For \unitcube{3} and \unitball{3} a similar construction works for cubic graphs $G$ where we
additionally require that the angles of edge-curves at the branching vertices are $120^\circ$.
The rest of the construction is analogous.\qed
\end{proof}

\section{Conclusions}

In this paper, we have determined the maximum cop number of circle, circular arc, function and
interval filament graphs, and we gave bounds for outer-string graphs, string graphs, and string
graphs on both orientable and non-orientable bounded genus surfaces. The following open problems remain.

\begin{problem}
Improve lower and upper bounds for the maximum cop number of string graphs, outer string graphs and
other intersection graphs of arc-connected sets in the plane such as 2-dimensional segments, boxes,
disks, unit disks, convex sets, etc.
\end{problem}

We note that it is proved in~\cite{unit_disk} that the maximum cop number of unit disk graphs is at
most 9. Their strategy is similar to our strategy for string graphs, by
applying Lemma~\ref{lem:guarding_shortest_path_neighbourhood}. The difference is that further geometric
properties of intersections of unit disks are proved which allows to guard a neighborhood of a
shortest path with just 3 cops instead of 5.

%
%
%

\bibliographystyle{plainnat}
\bibliography{cops_on_intersection_graphs}

\end{document}